\documentclass[12pt]{article}
\usepackage{amssymb,amsmath}
\usepackage{graphicx}
\usepackage{amsfonts}
\usepackage{float}
\usepackage{flafter}
\usepackage{color}

\labelwidth\leftmargini\advance\labelwidth-\labelsep

\def\R{\relax\ifmmode I\!\!R\else$I\!\!R$\fi}

\def\Z{\relax\ifmmode Z\!\!\!Z\else$Z\!\!\!Z$\fi}

\def\C{\relax\ifmmode C\!\!\!\!I\else$C\!\!\!\!I$\fi}

\def\K{\relax\ifmmode I\!\!K\else$I\!\!K$\fi}

\def\N{\relax\ifmmode I\!\!N\else$I\!\!N$\fi}

\newcounter{defcounter}[section]
{\vspace{0.1cm}\begin{sloppypar}\noindent\stepcounter{defcounter}{\bfseries
Definition
      \thesection.\thedefcounter}}%
{\end{sloppypar}\vspace{0.1cm}}

\newtheorem{theorem}{Theorem}[section]

\newtheorem{proposition}{Proposition}[section]

\newcommand{\proof}{{\bf Proof.} }

\newcommand{\qed}{\hfill $\square$}

\begin{document}
\thispagestyle{empty}
\begin{center}
{\Large {\bf Expansions of real numbers with respect to two integer bases}}
\end{center}
\begin{center}J\"org Neunh\"auserer\\
joerg.neunhaeuserer@web.de
\end{center}
\begin{center}
\begin{abstract}
We introduce and study expansions of real numbers with respect to two integer bases.
\end{abstract}
\end{center}
\section{Introduction}
We construct an expansion of real numbers in $[0,1]$ by blowing up the first basic interval $[0,1/b]$ in the classical $b$-adic expansion to $[0,1/a]$ with $a<b$. The details of the construction are given in the next section. Due to the overlaps of basic intervals the properties of these new expansions of real numbers turn out to be similar to the properties of the meanwhile classical $\beta$-expansions, $\sum d_{i}\beta^{-i}$ for $\beta\in(1,2)$. The literature on $\beta$-expansion is modestly growing since the pioneering work of Alfred Renyi (1921-1970) in 1950th. We refer to \cite{[SI2]} for an introduction to this topic. ~\\
In section 3 we will prove Sidorov's theorem \cite{[SI]} in the context of $a,b$-expansions using an ergodic theoretical approach: Almost all $x\in[0,1]$ have a continuum of such expansions. If we assume $b\ge a^2$, this is even true for all $x\in(0,1)$. The condition here is analogue of the condition $\beta<(\sqrt{5}+1)/2$, which leads to a continuum of $\beta$-expansion for all $x\in(0,1/(\beta-1))$. This has been proved by Paul Erd\"os (1913-1996) and his Hungarian collaborators, see \cite{[E2],[E3]}.\\
Now it is natural to ask, if there are exceptional $x\in(0,1)$ having an unique $a,b$-expansions if $a<b<a^2$. We will show in section 4 that there are at least countable many such exceptions, if $b<a^2$ and that there is an uncountable set with positive Hausdorff dimension of such exceptions if $b<a^2-2-(a\lceil b/a\rceil-b)$. The elementary number theoretical condition here is an analogue of the condition $\beta>\kappa(2)$, where $\kappa(2)$ is the transcendental Komornik-Loreti constant, which leads to a set of numbers with positive Hausdorff dimension, which have a unique $\beta$-expansion, see \cite{[VK]}. To characterize numbers which have a unique $a,b$-expansion, we again use ideas from dynamical systems.\\
In the last section we consider sets of numbers with digits from a proper subset of $\{0,\dots,b-1\}$ in their $a,b$-expansions. Hundred years ago Felix Hausdorff (1868-1942) \cite{[HA]} introduced his dimension to distinguish the size of such sets for $b$-adic expansions. We will determine the Hausdorff dimension of these sets with respect to $a,b$-expansions using a resent result of Mike Hochmann, see \cite{[HO]}.
\section{$a,b$-expansions}
Fix integers $a,b$ with $b>a>1$. Consider the space of digits $\Sigma=\{0,\dots, b-1\}^{\mathbb{N}}$. For $d=(d_{i})\in\Sigma$ and $n\in\mathbb{N}$ let $0_{n}(d)$ be the number of zeros in the sequence $(d_{1},\dots,d_{n})$.
We consider the map $\pi:\Sigma \to [0,1]$ given by
\[ \pi(d)=\sum_{i=1}^{\infty}d_{i}\left(\frac{1}{a}\right)^{0_{i}(d)}\left(\frac{1}{b}\right)^{i-0_{i}(d)}.\]
If $\pi(d)=x$ we call $d$ an $a,b$-expansion of $x$.
Several times we will use another description of the map $\pi$, which we now describe. For $j=0,\dots,b-1$ let $T_{j}:[0,1]\to [0,1]$ be contractions given by
\[ T_{0}(x)=\frac{x}{a}\quad\mbox{ and }\quad T_{j}(x)=\frac{x}{b}+\frac{j}{b}\quad\mbox{ for }\quad j=1,\dots,b-1.\]
By induction we have
\[ T_{d_{1}}\circ\dots \circ T_{d_{n}}(x)=\left(\frac{1}{a}\right)^{0_{n}(d)}\left(\frac{1}{b}\right)^{n-0_{n}(d)}x+\sum_{i=1}^{n}d_{i}\left(\frac{1}{a}\right)^{0_{i}(d)}\left(\frac{1}{b}\right)^{i-0_{i}(d)}.\]
Hence for all $d\in\Sigma$ and all $x\in I$
\[ \pi(d)=\lim_{n\to\infty}T_{d_{1}}\circ\dots \circ T_{d_{n}}(x).\]
Obviously we have
\[ [0,1]=\bigcup_{j=0}^{b-1}T_{i}([0,1]),\]
hence all $x\in[0,1]$ have at least one $a,b$-expansion. Such an expansion may be calculated using a greedy algorithm.
\begin{figure}
\vspace{0pt}\hspace{70pt}\scalebox{0.5}{\includegraphics{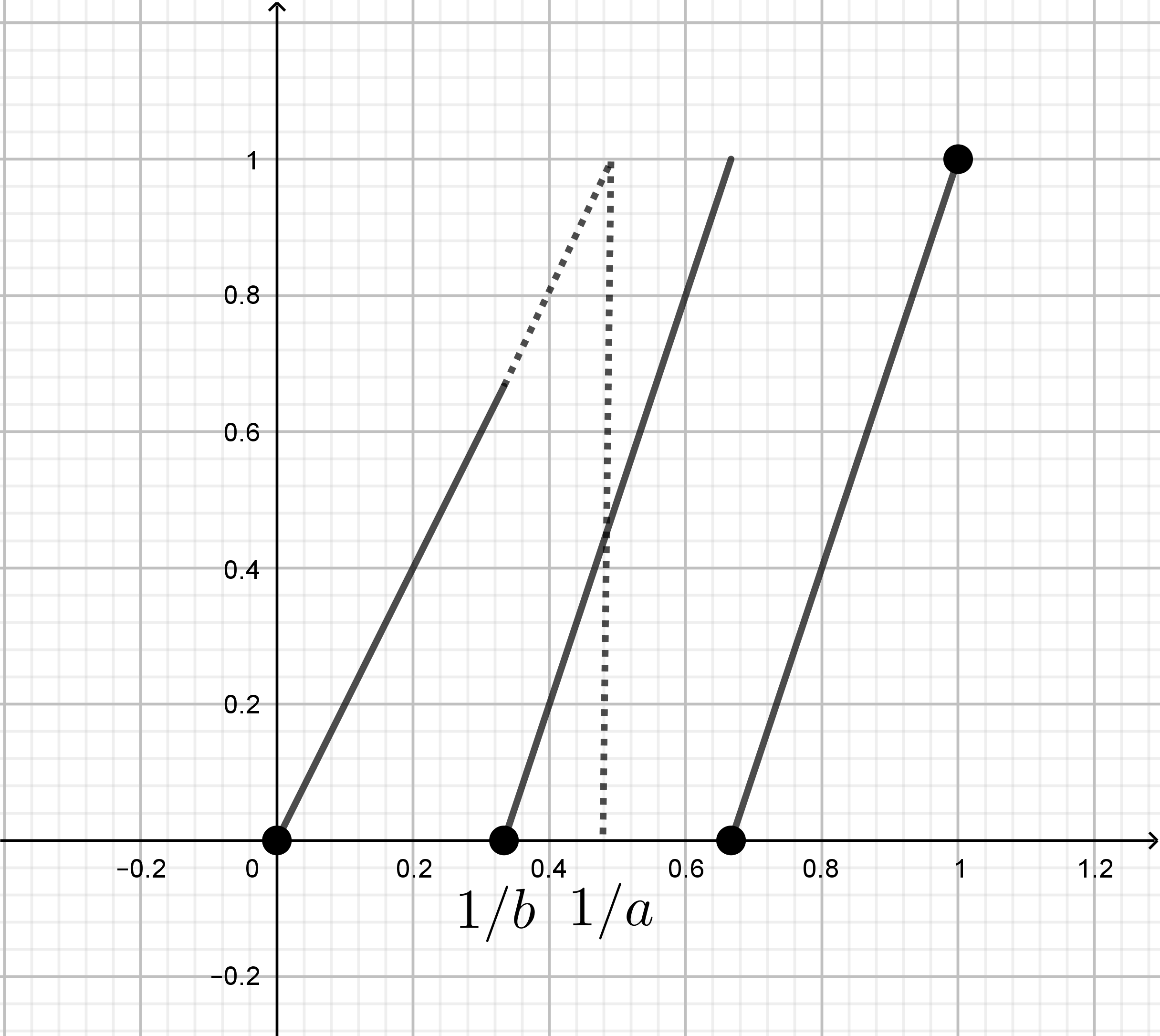}}
\caption{The map $G$ in the case $a=2,b=3$}
\end{figure}
Let $G:[0,1]\to [0,1]$ be given by
\[ G(x)=\left\{\begin{array}{ll} ax, & x\in [0,1/b) \\
         bx-j, & x \in[j/b,(j+1)/b)\quad \mbox{ for }\quad j=1,\dots,b-1 \\
         1, & \quad x=1 \end{array}\right.. \]
For $x\in [0,1]$ we define the greedy expansion $g=g(x)=(g_{i})\in\Sigma$ with respect to $a$ and $b$ by
\[ g_{i}=\lfloor b G^{i-1}(x)\rfloor,\]
where $\lfloor x\rfloor$ is the greatest integer not greater than $x$. We have
\begin{proposition}
For all $x\in [0,1]$ the greedy expansion $g(x)\in\Sigma$ with respect to $a$ and $b$  is an $a,b$-expansion of $x$, that is $\pi(g(x))=x$.
\end{proposition}
\proof
If $g_{i}=j$ we have $G^{i-1}(x)\in [j/b,(j+1)/b)$,  which implies $G^{i-1}(x)\in T_{g_{i}}([0,1])$. Since $G$ is defined by the inverse maps of $T_{i}$ we obtain
\[ x\in  T_{g_{1}}\circ\dots \circ T_{g_{i}}([0,1])\]
for all $i\in\mathbb{N}$, but this implies $\pi(g(x))=x$.
\qed
\section{A continuum of expansions}
Due to the intersection of the open intervals $(0,1/a)$ and $(1/b,2/b)$ in $(1/b,1/a)$ for $b>a$ the $a,b$-expansion of $x \in(0,1)$ will not be unique generically. We will first prove:
\begin{theorem}
Almost all $x\in(0,1)$ (in the sense of Lebesgue measure) have a continuum of $a,b$-expansions.
\end{theorem}
\proof
Let $G:[0,1]\to [0,1]$ be the map defined in the last section. $G$ is a piecewise linear expanding interval map and such maps are known to have an ergodic measure, which is equivalent to the Lebesgue measure, see \cite{[DE]} and \cite{[KH]}. By Poincare recurrence theorem for almost all $x\in [0,1]$ there is a $k\ge 0$ such that $G^{k}(x)\in (1/b,1/a)$. Note that each $x\in (1/b,1/a)$ has two different $a,b$-expansions, one starting with $0$ and other one starting with $1$.
 Hence for almost all $x\in [0,1]$ there is a $k\ge 0$ and a sequence $(d_{1},\dots,d_{k})\in\{0,\dots,b-1\}^{k}$ such that
 \[ x=T_{d_{1}}\circ \dots\circ T_{d_{k}}\circ T_{0}(x_{0})\quad\mbox{ and }\quad x=T_{d_{1}}\circ\dots\circ T_{d_{k}}\circ T_{1}(x_{1}),\]
where $x_{0},x_{1}\in (0,1)$ and $x_{0}\not=x_{1}$. For almost all $x$ both numbers $x_{1}(x),x_{2}(x)$ have two different $a,b$-expansion hence almost all $x$ have four different expansions. We use here that the intersection of two sets of full measure has full measure. Repeating this procedure $\aleph_{0}$ times we obtain $2^{\aleph_{0}}$ $a,b-$expansions for almost all $x\in [0,1]$, using the fact a countable intersection of sets of full measure has full measure.\qed~\\~\\
Under a suitable assumption on $a$ and $b$ we get a stronger result:
\begin{theorem}
If $a^2\le b$ all $x\in(0,1)$ have a continuum of $a,b$-expansions.
\end{theorem}
\proof Let $\Lambda_{0}=(\frac{1}{b},\frac{1}{a})$ and define recursively
\[ \Lambda_{n+1}=\Lambda_{n}\cup\bigcup_{j=0}^{b-1}T_{i}(\Lambda_{n}).\]
We prove
\[ \Lambda_{n}=\left(\frac{1}{a^{n}b},1+\frac{1}{a}(\frac{b-1}{b})^{n}-(\frac{b-1}{b})^{n}\right)\]
by induction. For $n=0$ the formula is obviously valid. Furthermore we have
\[ \Lambda_{n+1}=T_{0}(\Lambda_{n})\cup\Lambda_{n}\cup\bigcup_{j=1}^{b-1}T_{j}(\Lambda_{n})\]
\[=
 \left(\frac{1}{a^{n+1}b},\frac{1}{a}+\frac{1}{a^2}(\frac{b-1}{b})^{n}-\frac{1}{a}(\frac{b-1}{b})^{n}\right)\cup \left(\frac{1}{a^{n}b},1+\frac{1}{a}(\frac{b-1}{b})^{n}-(\frac{b-1}{b})^{n}\right)\]
\[\cup \left(\frac{b-1}{a^{n}b^2}+\frac{1}{b},1+\frac{1}{a}(\frac{b-1}{b})^{n+1}-(\frac{b-1}{b})^{n+1}\right)= \left(\frac{1}{a^{n+1}b},1+\frac{1}{a}(\frac{b-1}{b})^{n+1}-(\frac{b-1}{b})^{n+1}\right), \]
where we use $a^2\le b$ in the last equation to have an overlap of the intervals.\\
Now note that $\bigcup_{n\ge 0} \Lambda_{n}=(0,1)$. Hence for all  $x\in (0,1)$ there is a $k> 0$ and a sequence $(d_{1},\dots,d_{k})\in\{0,\dots,b-1\}^{k}$, such that
\[ x=T_{d_{1}}\circ\dots\circ T_{d_{k}}\circ T_{0}(x_{0})\mbox{ and }x=T_{d_{1}}\circ\dots\circ T_{d_{k}}\circ T_{1}(x_{1}),\]
where $x_{0},x_{1}\in (0,1)$ and $x_{0}\not=x_{1}$. Hence we obtain two expansions of $x$ that differ in the $k+1$ digit. Applying the result to $x_{0}(x)$ and $x_{1}(x)$ we obtain four expansions of $x$. Here we use that $x_{0}(x)$ and $x_{1}(x)$ are not at the boundary of $(0,1)$. Repeating this procedure $\aleph_{0}$ times we again see that there are $2^{\aleph_{0}}$ expansions of $x$.
 \qed~\\~\\
It is natural to ask if there are $x\in(0,1)$ that do not have a continuum of $a,b$-expansion if $a^2>b>a$. In the next section we will see that this is true.
\section{Unique expansions}
We consider the shift map $\sigma:\{0,\dots,b-1\}^{\mathbb{N}}\to \{0,\dots,b-1\}^{\mathbb{N}}$ given by $\sigma((d_{k}))=(d_{k+1})$. Using this map we characterise numbers which have a unique $a,b$-expansion as follows:
\begin{proposition}
The $a,b$-expansion $(d_{i})$ of $x\in[0,1]$ is unique if and only if \[ \pi(\sigma^{k}(d_{i}))\in ([0,1/b)\cup(1/a,1])\backslash\{j/b~|~j=2,\dots,b-1\}\]
for all $k\ge 0$.
\end{proposition}
$\pi(d)=\pi(e)$ with $d\not=e$ holds if and only if there exists a $k\ge 0$ such that $d_{k}\not=e_{k}$ and $\pi(\sigma^{k}((d_{i}))=\pi(\sigma^{k}((e_{i}))$. But this is equivalent to
\[\pi(\sigma^{k}(d_{i}))\in \bigcup_{l\not= m} T_{l}([0,1])\cap T_{m}([0,1])= [1/b,1/a]\cup\{j/b|j=2,\dots, b-1\}.\]
The proposition follows by contraposition.\qed~\\~\\
Obviously $0$ and $1$ have a unique $a,b$-expansion. Using the last proposition we are able to prove:
\begin{theorem}
If $a<b<a^2$ at least countable many $x\in(0,1)$ have a unique $a,b$-expansion. If $a<b<a^2-2-r$, where $r\in\{0,\dots,b-1\}$ is given by $b+r\in\{la|l\in\mathbb{N}\}$, uncountable many $x\in(0,1)$ have a unique $a,b$-expansion. Furthermore the Hausdorff dimension of the set of these numbers is positive.
\end{theorem}
\proof Let $d$ be the periodic sequence in $\Sigma$ given by $d=(d_{i})=(0,a,0,a,\dots)$. Using the geometric series we have
\[ \pi(d)=\frac{1}{1+a^{-1}b^{-1}}\cdot\frac{1}{b}<\frac{1}{b}.\]
On the other hand
\[ \pi(\sigma(d))= \frac{1}{1+a^{-1}b^{-1}}\cdot \frac{a}{b}>\frac{1}{a}\]
 by the assumption $a^2>b>a$ for integers $a,b$. Furthermore $\pi(\sigma(d))\not=j/b$ for $j=2,\dots,b-1$. By proposition 4.1 $x=\pi(d)\in(0,1)$ has an unique $a,b$-expansion. Considering sequences $d=(0,\dots,0,a,0,a,\dots)$ we see that there are at least countable many such $x\in(0,1)$.\\
Let $b=al-r$ with $r\in\{0,1\dots ,b-1\}$ with $l\in\mathbb{N}$ minimal and consider the set $V=\{0l,l0\}^{\mathbb{N}}$ in $\Sigma$. The set is not invariant under the shift map $\sigma$ but
\[ U=\bigcup_{k=0}^{\infty}\sigma^{k}(V)=V\cup(\{0\}\times V) \cup(\{l\}\times V)\]
is a $\sigma$-invariant set. The sequence $d\in U$ with $d_{1}=0$, that has the largest projection under $\pi$, obviously is $d=(0,l,0,0,l,0,l,0,\dots)$. We have
\[ \pi(d)=\frac{l}{ b}\left(\frac{1}{a}+\frac{1}{1+ab}\right)=\frac{b+r}{a^2b}+\frac{b+r}{b(1+ab)}\]
\[ =\frac{1}{a^2b}(b+r+\frac{a(b+r)}{1+ab})<\frac{1}{b}  \]
by our assumption $b+r+2<a^2$ and by the inequality $a(b+r)/(1+ab)<2$. The sequence $d\in U$ with $d_{1}=l$, that has the smallest projection under $\pi$, obviously is $d=(l,0,0,l,0,l,0,\dots)$. We have
\[ \pi(d)=\frac{l}{b}+\frac{l}{ab+a^2b^2}\ge \frac{1}{a}+\frac{1}{a^2+a^3b}> \frac{1}{a}  \]
since $l \ge b/a$. We have thus shown that $\pi(\sigma^{k}(U))=\pi(U)\subseteq [0,1/b)\cup(1/a,1]$ for all $k\ge 0$. We now prove that $\pi(U)$ is uncountable with positive Hausdorff dimension. The subset $\pi(V)\subseteq \pi(U)$ is the attractor of the iterated function system $(T_{0}\circ T_{l},T_{l}\circ T_{0})$. This means $T_{0}\circ T_{l}(\pi(V))\cup T_{l}\circ T_{0}((\pi(V))=\pi(V)$. The iterated functions system fulfills the open set condition, since $T_{0}\circ T_{l}((0,1))\cap T_{l}\circ T_{0}((0,1)))=\emptyset$. By \cite{[MO]} the set $\pi(V)$ (and hence $\pi(U)$ as well) has positive Hausdorff dimension and is thus uncountable. Note that $\pi(\sigma^{k}(d))=j/b$ for some $j$ and $k$ hold for sequences $d$ in a at most countable set $S\subseteq U$. Each $x\in\pi(U\backslash S)$ has unique $a,b$-expansion by proposition 4.1 and this set has positive Hausdorff dimension and is uncountable.
\qed~\\~\\
The upper bound in the second part of our theorem is due to the idea of the proof. We do not know if this bound is sharp.
\section{Restricted digits}
Let $D\subseteq\{0,\dots,b-1\}$ be subset of digits. We are interested in the set of numbers which have an $a,b$-expansion with digits in $D$, that is the set $\pi(D^{\mathbb{N}})$. We denote the  cardinality of $D$ by $|D|$. If $|D|=1$ the set $\pi(D^{\mathbb{N}})$ obviously contains only one point and we already know that $\pi(D^{\mathbb{N}})=[0,1]$ if $|D|=b$. Moreover we obtain:
\begin{theorem}
Let $D\subseteq\{0,\dots,b-1\}$ with $1<|D|<b$.
If $0\not\in D$ we have \[\dim_{H}\pi(D^{\mathbb{N}})=\frac{\log|D|}{\log b}.\]
If $0\in D$ and $\log(b)/\log(a)\not\in\mathbb{Q}$ we have $\dim_{H}\pi(D^{\mathbb{N}})=\min\{s,1\}$ where $s$ is given by
\[ \left(\frac{1}{a}\right)^{s}+(|D|-1)\left(\frac{1}{b}\right)^{s}=1.\]
The same is true if $\min(D\backslash\{0\})\ge b/a$ no matter where $\log(b)/\log(a)\in\mathbb{Q}$.
\end{theorem}
\proof If $0\not\in D$ the set $\pi(D^{\mathbb{N}})$ contains the numbers in $[0,1]$, which have a $b$-adic expansion with digits in $D$. The dimension formula goes back to Hausdorff \cite{[HA]} in this case. Now assume $0\in D$ and $\log(b)/\log(a)\not\in\mathbb{Q}$. The assumption means that $b^{n}\not=a^{m}$ for all $n,m\ge 1$. The identity of the maps
\[  T_{d_{1}}\circ\dots \circ T_{d_{n}}=T_{g_{1}}\circ\dots \circ T_{g_{m}}\]
for some sequences $d=(d_{1},\dots, d_{n})$ and $g=(g_{1},\dots, g_{m})$ is equivalent to
\[\left(\frac{1}{a}\right)^{0_{n}(d)}\left(\frac{1}{b}\right)^{n-0_{n}(d)}x+\sum_{i=1}^{n}d_{i}\left(\frac{1}{a}\right)^{0_{i}(d)}\left(\frac{1}{b}\right)^{i-0_{i}(d)}\]
\[=\left(\frac{1}{a}\right)^{0_{m}(g)}\left(\frac{1}{b}\right)^{m-0_{m}(g)}x+\sum_{i=1}^{m}g_{i}\left(\frac{1}{a}\right)^{0_{i}(g)}\left(\frac{1}{b}\right)^{i-0_{i}(g)}\]
for all $x\in[0,1]$. By our assumption this implies $n=m$ and $0_{n}(d)=0_{m}(g)$. Furthermore we obtain $d_{i}=g_{i}$ for $i=1,\dots ,n$. This means that there are no exact overlaps of intervals of the form $T_{d_{1}}\circ\dots\circ T_{d_{n}}([0,1])$. By Hochmann's result  \cite{[HO]} this implies $\dim_{H}\pi(D^{\mathbb{N}})=\min\{s,1\}$, where $s$ is the similarity dimension of the iterated function system $\{T_{i}|i\in D\}$. That means
\[ \sum_{i\in D}|T_{i}^{'}|^{s}=\left(\frac{1}{a}\right)^{s}+(|D|-1)\left(\frac{1}{b}\right)^{s}=1.\]
If $\min(D\backslash\{0\})\ge b/a$, we have $T_{i}((0,1))\cap T_{j}((0,1))=\emptyset$ for $i\not=j$. In this case the result follows from the classical theory of iterated function systems fulfilling the open set condition, see \cite{[MO]}, or again from \cite{[HO]}.
\qed~\\~\\
If we have $\log(b)/\log(a)\in\mathbb{Q}$ and $\min(D\backslash\{0\})< b/a$, there may occur exact overlaps of intervals $T_{d_{1}}\circ\dots\circ T_{d_{n}}([0,1])$. In this situation a combinatorial approach is necessary to calculate the Hausdorff dimension of $\pi(D^{\mathbb{N}})$. We refer here to \cite{[NW]}.\\


\begin{thebibliography}
\small
\bibitem{[DE]} M. Denker, Introduction to analysis of dynamical systems, Springer, Berlin, 2005.
\bibitem{[E2]} P. Erd\"os, I. Joo, On the number of expansions $\sum_{i=1}^{\infty} q^{-n_{i}}$, Ann. Univ. Sci.
Budapest 35 (1992), 129-132.
\bibitem{[E3]} P. Erd\"os, I. Joo and V. Komornik, On the number of $q$-expansions, Ann. Univ. Sci. Budapest
E\"otv\"os Sect. Math. 37 (1994), 109-118.
\bibitem{[HA]} F. Hausdorff, Dimension und \"ausseres Maß. Math. Annalen 79 (1919), 157-179.
\bibitem{[HO]} M. Hochmann, On self-similar sets with overlaps and inverse theorems for entropy, Annals of Mathematice 180 (2014), no. 2, 773-822.
\bibitem{[KH]} A. Katok and B. Hasselblatt, Introduction to the Modern Theory of Dynamical Systems, Cambridge University Press, 1995.
\bibitem{[MO]} P. Moran, Additive functions of intervals and Hausdorff measure, Proc. Cambridge
Philos. Soc. 42 (1946) 15-23.
\bibitem{[NW]} S.M. Ngai and Y. Wang, Journal of the London Mathematical Society, vol. 63 (2001), 655-672.
\bibitem{[RE]} A. Renyi, Representations for real numbers and their ergodic properties, Acta Math.
Acad. Sci. Hung. 8 (1957) 477-493.
\bibitem{[VK]} de Vries M. and Komornik V., Unique expansions of real numbers. Adv. Math. 221 (2009), 390-427.
\bibitem{[SI]} N. Sidorov, Almost every number has a continuum of $\beta$-expansions, Amer.Math. Monthly 110 (2003), 838-842.
\bibitem{[SI2]} N. Sidorov, Arithmetic dynamics, Lond. Math. Soc. Lect. Note Ser. 310 (2003), 145-189.
\end{thebibliography}
\end{document}